\documentclass[10pt]{amsart}
\usepackage[square,compress,comma, numbers]{natbib}
\usepackage[colorlinks=true, citecolor=blue, linkcolor=blue]{hyperref}
\usepackage{amsfonts}
\allowdisplaybreaks[4]
\usepackage{amssymb}
\usepackage{amsmath}
\usepackage{color}
\newcommand{\abs}[1]{\left\lvert #1 \right\rvert}

\def\E#1{\mathbb{E}\left \{#1 \right\}}

\definecolor{c20}{rgb}{0.,0.7,0.}
\definecolor{c30}{rgb}{0.,0.,1.}
\definecolor{c40}{rgb}{1,0.1,0.7}
\definecolor{c50}{rgb}{1,0,0}
\definecolor{c60}{rgb}{1,0.9,0.1}
\definecolor{c70}{rgb}{0.50,1.00,0.00}

\def\N{\mathbb{N}}

\numberwithin{equation}{section}
\newtheorem{theo}{Theorem}[section]
\newtheorem{sat}[theo]{Proposition}
\newtheorem{de}[theo]{Definition}
\newtheorem{lem}{Lemma}[section]

\newtheorem{korr}[theo]{Corollary}
\newtheorem{remark}[theo]{Remark}
\newtheorem{remarks}[theo]{Remarks}
\numberwithin{equation}{section}

\newcommand{\prooftheo}[1]{ \textsc{Proof of Theorem} \ref{#1} }

\newcommand{\pk}[1]{\mathbb{P} \left\{ #1 \right\} }

\newcommand{\QED}{\hfill $\Box$}
\newcommand{\COM}[1]{}
\def\IF{\infty}
\newcommand{\R}{\mathbb{R}}

\newcommand{\BQN}{\begin{eqnarray}}
\newcommand{\EQN}{\end{eqnarray}}
\newcommand{\BQNY}{\begin{eqnarray*}}
\newcommand{\EQNY}{\end{eqnarray*}}

\def\polhk#1{\setbox0=\hbox{#1}{\ooalign{\hidewidth
\lower1.5ex\hbox{`}\hidewidth\crcr\unhbox0}}}

\def\cL#1{\textcolor{c50}{#1}}
\def\cL#1{#1}

\newcommand{\netheo}[1]{{Theorem \ref{#1}}}

\def\rw{\rightarrow}

\def\IF{\infty}

\date{}

\def\LT{\left}
\def\RT{\right}

\def\rw{\rightarrow}

\def\TT{\mathcal{T}}

\def\vn{\varepsilon}
\def\Var{\text{Var}}

\def\PP{\widetilde{\mathcal{P}}}

\newcommand{\limit}[1]{\lim_{#1 \to \infty}}

\newcommand{\BS}{\begin{sat}}
\newcommand{\ES}{\end{sat}}
\newcommand{\BT}{\begin{theo}}
\newcommand{\ET}{\end{theo}}
\newcommand{\BK}{\begin{korr}}
\newcommand{\EK}{\end{korr}}

\newcommand{\BD}{\begin{de}}
\newcommand{\ED}{\end{de}}
\newcommand{\BIT}{\begin{itemize}}
\newcommand{\EIT}{\end{itemize}}
\newcommand{\BDI}{\begin{description}}
\newcommand{\EDI}{\end{description}}
\newcommand{\BRM}{\begin{remarks}}
\newcommand{\ERM}{\end{remarks}}
\newcommand{\BEL}{\begin{lem}}
\newcommand{\EEL}{\end{lem}}

\def\TT{\mathcal{T}}
\def\TT{\R }

\def\rw{\rightarrow}
\def\LT{\left}
\def\RT{\right}

\def\Var{\text{Var}}
\def\vn{\varepsilon}

\def\LT{\left}
\def\RT{\right}

\def\rw{\rightarrow}

\def\TT{\mathcal{T}}

\def\vn{\varepsilon}
\def\Var{\text{Var}}

\def\LT{\left}
\def\RT{\right}

\def\rw{\rightarrow}

\def\TT{\mathcal{T}}

\def\vn{\varepsilon}
\def\Var{\text{Var}}

\def\PP{\widetilde{\mathcal{P}}}

\def\cL#1{\textcolor{c50}{#1}}
\def\cL#1{#1}

\begin{document}
\title{Parisian Ruin of Brownian Motion Risk Model over an Infinite-time Horizon}
\author{Long Bai}
\address{Long Bai, Department of Actuarial Science, University of Lausanne, UNIL-Dorigny, 1015 Lausanne, Switzerland
}
\email{Long.Bai@unil.ch}

\bigskip
\maketitle
{\bf Abstract:}
Let $B(t), t\in \mathbb{R}$ be a standard Brownian motion. In this paper, we derive the exact asymptotics of the probability of Parisian ruin on infinite time horizon for the following risk process
\begin{align}\label{Rudef}
R_u^{\delta}(t)=e^{\delta t}\left(u+c\int^{t}_{0}e^{-\delta v}d v-\sigma\int_{0}^{t}e^{-\delta v}d B(v)\right),\quad  t\geq0,
\end{align}
 where $u\geq 0$ is the initial reserve, $\delta\geq0$ is the force of interest, $c>0$ is the rate of premium and $\sigma>0$ is a volatility factor.
 Further, we show the asymptotics of the Parisian ruin time of this risk process.

\def\TTT{\mathcal{T}}
\def\TT{\mathcal{T}}
\def\Z{\mathbb{Z}}
\def\HWD{\mathcal{H}_\alpha^\delta}
\def\phd{\mathcal{P}_{\alpha,\delta}^h}
\section{Introduction}
In the risk theory, the surplus process of an insurance company  can be modeled by
\BQNY
R_u(t)=u+ct-X(t), \quad t\ge0,
\EQNY
see \cite{MR1458613}, where $u\geq 0$ is the initial reserve, $ct$ models the total premium received up to time $t$, and $ X(t),t\geq0 $ denotes the aggregate claims process. In \cite{dkebicki2015parisian,debicki2015parisian}, the Parisian ruin of $R_u(t)$ is defined by
\BQN\label{nintere}
\mathcal{P}_S(u, T_u)=\pk{\inf_{t\in[0,S]} \sup_{s\in[t,t+T_u]} R_u(s)<0}, \quad S\in(0,\IF],
\EQN
where $T_u$ models the pre-specified time which is a function of $u$. For $X(t), t\geq 0$ a Gaussian process, the asymptotics of $\mathcal{P}_S(u, T_u)$ over finite-time horizon, i.e. $S\in(0,\IF)$, is investigated in \cite{debicki2015parisian}. Further, \cite{dkebicki2015parisian} showed the tail asymptotic results of $R_u(t)$ over infinite-time horizon, i.e. $S=\IF$ in \eqref{nintere}, where $X(t)$ is a  self-similar Gaussian process.
In this paper considering the nature of the financial market,  we introduce the force of interest $\delta$ into the model $R_u(t)$ as $R^\delta_u(t)$ in \eqref{Rudef} when $X(t)=B(t)$. \cite{ParisianBrownianfinite2017} gave an approximation of the Parisian ruin probability
\begin{align*}
\mathcal{K}_S^{\delta}(u,T_u):=\pk{\inf_{t\in[0,S]} \sup_{s\in[t,t+T_u]} R_u^{\delta}(s)<0}, \quad S\in(0,\IF),
\end{align*}
as $u\rightarrow\infty$. See \cite{rolski2009stochastic, DHJ13a,HX2007} for more studies on risk models with force of interest.
In the literature, no results are available for the approximation of Parisian ruin probability over infinite time horizon for $\delta>0$.
%But the Parisian ruin probabilities for the models with force of interest over infinite time horizon have never  been covered by the aforementioned literatures.
 In this contribution we shall investigate the asymptotics of the Parisian ruin probability
\begin{align*}
\mathcal{K}^{\delta}(u,T_u):=\mathbb{P}\left\{\inf_{t\geq 0} \sup_{s\in[t,t+T_u]} R_u^{\delta}(s)<0\right\},
\end{align*}
as $u\rightarrow\infty$ where $T_u\geq 0$  models the pre-specified time satisfying
\BQN\label{TT}
\limit{u} T_u= T\in[0,\IF]. 
\EQN
 When $\delta=0$ and $T\in[0,\IF)$, \cite{dkebicki2015parisian} showed that (hereafter $\sim$ means asymptotic equivalence)
 \begin{align*}
\mathcal{K}^0(u,T_u)=\mathbb{P}\left\{\inf_{t\geq 0} \sup_{s\in[t,t+T_u]} \LT(u+cs-\sigma B(s)\RT)<0\right\}
\sim\mathcal{F}\LT(\frac{2c^2T}{\sigma^2}\RT)\exp\LT(-\frac{2cu}{\sigma^2}\RT),
\ u\rw\IF,
\end{align*}
where
\BQNY
\mathcal{F}\LT(T\RT)=\lim_{\lambda\rw\IF}\frac{1}{\lambda}\E{\sup_{t\in[0,\lambda]}\inf_{ s\in[0,T]} e^{ \sqrt{2}B(t+s)-(t+s)}}.
\EQNY
Hereafter we make the convention that $\sup{\{\emptyset\}}=0$ and $\inf{\{\emptyset\}}=\IF$.\\
Complementary, we investigate the conditional distribution of the ruin time for the surplus process $R_u^{\delta}(t)$. The classical ruin time, e.g., \cite{DHJ13a, HJ13, MR2462285}, is defined as
\BQN\label{fo21}
\tau(u)=\inf\{t>0: R_u^{\delta}(t)<0\}.
\EQN
Here as in \cite{dkebicki2015parisian,ParisianBrownianfinite2017} we define the Parisian ruin time of  the  risk process $R_u^{\delta}(t)$  by %in \cite{czarna2011ruin, loeffen2013parisian, dkebicki2015parisian, debicki2015parisian} by
\BQN\label{eq:eta2}
\eta(u)=\inf\{t\ge T_u: t-\kappa_{t,u}\ge T_u, R_u^{\delta}(t)<0\},\ \ \ \text{with} \ \kappa_{t,u}=\sup\{s\in[0,t]: R_u^{\delta}(s)\ge0\},
\EQN
and $\tau(u)=\eta(u)$ when $T_u\equiv 0$.\\
Brief outline of the rest of the paper: In Section 2 we present our main results on the asymptotics of $\mathcal{K}^{\delta}(u,T_u)$ as $u\rw\IF$ and the approximation of the Parisian ruin time. All the  proofs are relegated to Section 3.

\section{Main results}
Before giving the main results, we shall introduce a constant as
\BQN\label{eq:GP}
\PP^f_a[0,\IF)=\lim_{\lambda\to\IF} \PP^f_a[0,\lambda]\in(0,\IF),
\EQN
with
\BQNY
\PP^f_a[0,\lambda]=\E{\sup_{t\in[0,\lambda]}\inf_{ s\in[a,1]} \exp\LT( \sqrt{2}B(st)-st-f(st)\RT)}\in(0,\IF),
\EQNY
 where $\lambda\ge0, a\in[0,1]$ and $f(t)$ is a  continuous function satisfying $\lim_{t\rw\IF}\frac{f(t)}{t^\epsilon}=\IF$ for some $\epsilon>0$.

Note further that $\PP^f_0[0,\lambda]=e^{-f(0)}$ and $$\PP^f_1[0,\lambda]=\E{\sup_{t\in[0,\lambda]} \exp\LT( \sqrt{2}B(t)-\abs{t}-f(t)\RT)},$$
see e.g. \cite{MR1993262,GeneralPit16,MR3317852} for the bounds of $\PP^f_a[0,\IF)$ and more details. \\
\COM{It is known that$\mathcal{H}_{1}=1$ , $\mathcal{H}_{2}=1/\sqrt{\pi}$ ,$\mathcal{P}_{1}^{1}=1+\frac{1}{a}$ and $\mathcal{P}_{2}^{1}=\frac{1}{2}\left(1+\sqrt{1+\frac{1}{a}}\right)$, see \cite{MR1993262,GeneralPit16,MR3317852}.}
\COM{Through this paper $\sim$ means asymptotic equivalence when the argument tends to $0$ or $\IF$.}
Recall that $\Phi(\cdot),\Psi(\cdot)$ denote the distribution function and the survival function of an $\mathcal{N}(0,1)$ random variable, respectively, and  $\Psi(u)\sim\frac{1}{\sqrt{2\pi}u}e^{-\frac{u^2}{2}}, u\rw\IF$.
\COM{\BT\label{pan1}
For any $T\in(0,\infty)$ and $\delta>0$, we have
\BQN
\psi_{T}(u)\sim
\frac{\sigma}{\sqrt{\delta\pi}}\frac{\sqrt{1-e^{-2\delta T}}}{u+\frac{c}{\delta}(1-e^{-\delta T})}\exp\left(-\frac{\delta(u+\frac{c}{\delta}(1-e^{-\delta T}))^2}{\sigma^2(1-e^{-2\delta T})}\right),
\EQN
as $u\rw\IF$.
\ET}
\BT\label{paripan1}
For $ \delta>0$ and $ T_u $ satisfying \eqref{TT}, we have
\begin{align*}
\mathcal{K}^{\delta}(u,T_u)&\sim\E{\sup_{t\in [0, \IF)}\inf_{s\in[a,1]}
\exp\LT(\sqrt{2}B(st)-st
-\LT(\sqrt{st}-\frac{c}{\sigma\sqrt{\delta}}\RT)^2\RT)}
\Psi\LT(\frac{1}{\sigma}\sqrt{2\delta u^2+4cu}\RT)\\
&=\PP^f_a[0,\IF)\Psi\LT(\frac{1}{\sigma}\sqrt{2\delta u^2+4cu}\RT),\ u\rw\IF,
\end{align*}
where $a=e^{-2\delta T}$ and $f(t)=\LT(\sqrt{t}-\frac{c}{\sigma\sqrt{\delta}}\RT)^2$.
\ET
\begin{remark}
In \netheo{paripan1}, if $T=0$, $a=1$,we get the asymptotic result of the classical ruin probability, i.e., as $u\rw\IF$
\BQNY
\mathcal{K}^{\delta}(u,0)=\mathbb{P}\left\{\inf_{t\geq 0} R_u^{\delta}(s)<0\right\}
\sim\E{\sup_{t\in [0, \IF)}
\exp\LT(\sqrt{2}B(t)-t
-\LT(\sqrt{t}-\frac{c}{\sigma\sqrt{\delta}}\RT)^2\RT)}
\Psi\LT(\frac{1}{\sigma}\sqrt{2\delta u^2+4cu}\RT)
\EQNY
which corresponds to the results in \cite{Threshold2016}.\\
Moreover, according to \cite{RPwCA1977} (see also \cite{ADARP1975})  we have
\BQN\label{CC}
\mathcal{K}^{\delta}(u,0)=\Psi\LT(\frac{\sqrt{2\delta}}{\sigma}\LT(u+\frac{c}{\delta}\RT)\RT)\Big/\Psi\LT(\frac{\sqrt{2}c}{\sigma\sqrt{\delta}}\RT).
\EQN
\end{remark}
\BT\label{ruintime}
  Let $\eta(u)$ satisfy (\ref{eq:eta2}), under the assumptions and notation of \netheo{paripan1}, we have for $\delta>0$ and $x\in\LT(-\frac{c^2}{\delta^2},\IF\RT)$
\BQN\label{parifo4}
\pk{u^2\LT(e^{-2\delta\eta(u)}-\LT(\frac{c}{\delta u+c}\RT)^2\RT)\leq x\ \big\lvert \eta(u)<\IF}\sim\frac{\PP^f_a[0,\frac{c^2}{\sigma^2\delta}+\frac{\delta x}{\sigma^2}]}{\PP^f_a[0,\IF)},
\ u\rightarrow\IF.
\EQN
\ET
\begin{remarks}
i) When $\delta=0$, \cite{dkebicki2015parisian} showed that for $x\in\R$
\BQNY
\pk{u^{-\frac{1}{2}}\LT(\eta(u)-\frac{u}{c}\RT)\leq x\ \big\lvert \eta(u)<\IF}\sim\Phi(c x), \ u\rightarrow\IF.
\EQNY
ii) When $T_u\equiv0$, $\eta(u)=\tau(u)$, by \eqref{parifo4}, we have
\BQNY
\pk{u^2\LT(e^{-2\delta\tau(u)}-\LT(\frac{c}{\delta u+c}\RT)^2\RT)\leq x\ \big\lvert \eta(u)<\IF}\sim\frac{\PP^f_1[0,\frac{c^2}{\sigma^2\delta}+\frac{\delta x}{\sigma^2}]}{\PP^f_1[0,\IF)},
\ u\rightarrow\IF,
\EQNY
which corresponds to the result in \cite{Threshold2016}.
\end{remarks}
\section{Proofs}
Hereafter we assume that $ \mathbb{C}_i, i\in \N$ are positive constants.

\prooftheo{paripan1}
We have for $u>0$
\BQNY
\mathcal{K}^{\delta}(u,T_u)=\pk{\inf_{t\in[0,\IF)}\sup_{s\in[t,t+T_u]}R_u^{\delta}(s)<0}=\pk{\inf_{t\in[0,\IF)}\sup_{s\in[t,t+T_u]}\widetilde{R}_u^{\delta}(s)<0},
\EQNY
where
$$\widetilde{R}_u^{\delta}(s)=u+c \int_0^s e^{-\delta v}dv-\sigma \int_0^s e^{-\delta v}dB(v),\ \ t\ge0.$$
Since for $t\in(0,\IF)$
$$\E{\LT[\sigma\int_0^t e^{-\delta v}dB(v)\RT]^2}=\frac{\sigma^2}{2\delta}\LT(1-e^{-2\delta t}\RT),$$
then
$$\sup_{t\in[0,\IF)}\E{\LT[\sigma\int_0^t e^{-\delta v}dB(v)\RT]^2}<\IF$$
implies that $$ \sup_{t\in[0,\IF)}\E{\LT| \sigma\int_0^t e^{-\delta v}dB(v)\RT|}<\IF,$$
by the martingale convergence theorem, see \cite{PP1966}, $\widetilde{R}_u^{\delta}(\IF):=\lim_{t\rw\IF}\widetilde{R}_u^{\delta}(t)$ exists and is finite almost surely. Thus for any $u>0$
\begin{align*}
\psi(u):=&\pk{\inf_{t\in[0,\IF)}\sup_{s\in[t,t+T_u]}\widetilde{R}_u^{\delta}(s)<0}=\pk{\inf_{t\in[0,\IF]}\sup_{s\in[t,t+T_u]}\widetilde{R}_u^{\delta}(s)<0}\\
=&\pk{\sup_{ t\in[0,\IF]}\inf_{s\in[t,t+T_u]}\LT(\sigma\int_0^s e^{-\delta v}dB(v)-c\int_0^s e^{-\delta v}dv\RT)>u}.
\end{align*}
Using a change of variable $s=-\frac{1}{2\delta}\ln s^*, s^*\in[t^*e^{-2\delta T_u},t^*],\ t^*\in[0,1]$, we have
\begin{align*}
\psi(u)&=\mathbb{P}\left\{\sup_{t^*\in[0,1]}\inf_{s^*\in[t^*e^{-2\delta T_u},t^*]}\LT(\sigma\int_0^{-\frac{1}{2\delta}\ln s^*} e^{-\delta v}dB(v)-c\int_0^{-\frac{1}{2\delta}\ln s^*} e^{-\delta v}dv\RT)>u\right\}\\
&=\mathbb{P}\left\{\sup_{t^*\in[0,1]}\inf_{s^*\in[t^*e^{-2\delta T_u},t^*]}\LT(\sigma\int_0^{-\frac{1}{2\delta}\ln s^*} e^{-\delta v}dB(v)-\frac{c}{\delta}(1-{s^*}^{\frac{1}{2}})\RT)>u\right\}.
\end{align*}
For simplicity, we still use $s,t$ instead of $s^*, t^*$.

Below, we set $Z(s)=\sigma\int_0^{-\frac{1}{2\delta}\ln s} e^{-\delta v}dB(v)$
with variance function given by
\BQNY
V_Z^2(s)=Var\left(\sigma\int_0^{-\frac{1}{2\delta}\ln s} e^{-\delta v}d B(v)\right)=\frac{\sigma^2}{2\delta}(1-s), \quad s\in [0,1].
\EQNY
We show next that for $u$ sufficiently large
$$M_u(t):=\frac{uV_Z(t)}{G_u(t)}=\frac{\frac{\sigma}{\sqrt{2\delta}}\sqrt{1-t}}{1+\frac{c}{\delta u}(1-t^{1/2})},  \quad 0\leq t \leq 1,$$
with $G_u(t):=u+\frac{c}{\delta}(1-t^{\frac{1}{2}})$ attains its maximum at the unique point
$$t_u=\LT(\frac{c}{\delta u+c}\RT)^2.
$$
%We calculate the derivative of $M_u(t)$, i.e,
In fact, we have for $t\in(0,1)$
\BQN\label{VIF1}
[M_u(t)]_t:=\frac{d M_u(t)}{dt}&=&\frac{dV_Z(t)}{dt}\cdot\frac{u}{G_u(t)}-\frac{V_Z(t)}{G_u^2(t)}\LT(-\frac{cu}{2\delta}t^{-\frac{1}{2}}\RT)\nonumber\\
&=&\frac{u}{2G_u^2(t) V_z(t)}\left[\frac{dV_Z^2(t)}{dt}G_u(t)+V_Z^2(t)\frac{ct^{-\frac{1}{2}}}{\delta}\right]\nonumber\\
&=&\frac{u\sigma^2 t^{-1/2}}{4\delta G_u^2(t) V_Z(t)}\LT[\frac{c}{\delta}-\LT(u+\frac{c}{\delta}\RT)t^{\frac{1}{2}}\RT].
\EQN
Letting $[M_u(t)]_t=0$, we get $t_u=\LT(\frac{c}{\delta u+c}\RT)^2$.\\
By (\ref{VIF1}), $[M_u(t)]_t>0$ for $t\in(0,t_u)$ and $[M_u(t)]_t<0$ for $t\in(t_u,1)$, so $t_u$ is the unique maximum point of $M_u(t)$ over $[0,1]$. Further
\BQNY
M_u:=M_u(t_u)=\frac{\sigma u}{\sqrt{2\delta u^2+4cu}}=\frac{\sigma}{\sqrt{2\delta}}(1+o(1)),\ u\rw\IF.
\EQNY
Set $\delta(u)=\LT(\frac{\ln u}{u}\RT)^2$,
$\Delta(u)=[0,t_u+\delta(u)]$ and for some positive constant $\lambda$
\begin{align*}
I_u(k)=\LT[k\lambda u^{-2},(k+1)\lambda u^{-2}\RT],\  \cL{k\in \N},\ \
 N(u)=\LT\lfloor \lambda^{-1}(\ln u)^{2}\RT\rfloor.
\end{align*}
We have for $u$ large enough
\BQN
&&\psi(u)\geq\pk{\sup_{t\in [0, t_u+\lambda u^{-2}]}\inf_{s\in[te^{-2\delta T_u},t]}\LT(\sigma\int_0^{-\frac{1}{2\delta}\ln s} e^{-\delta v}dB(v)-\frac{c}{\delta}(1-{s}^{\frac{1}{2}})\RT)>u}=:\Pi_0(u),\label{bound1}\\
&&\psi(u)\leq \Pi_0(u)+\Pi_1(u)+\Pi_2(u)+\Pi_3(u),\label{bound2}
\EQN
where for $\theta\in(0,1)$
\begin{align*}
&\Pi_1(u)=\sum_{i=1}^{N(u)}\pk{\sup_{t\in (t_u+I_u(k))}\inf_{s\in[te^{-2\delta T_u},t]}\LT(\sigma\int_0^{-\frac{1}{2\delta}\ln s} e^{-\delta v}dB(v)-\frac{c}{\delta}(1-{s}^{\frac{1}{2}})\RT)>u},\\
&\Pi_2(u)=\pk{\sup_{t\in [t_u+\delta(u),\theta]}\inf_{s\in[te^{-2\delta T_u},t]}\LT(\sigma\int_0^{-\frac{1}{2\delta}\ln s} e^{-\delta v}dB(v)-\frac{c}{\delta}(1-{s}^{\frac{1}{2}})\RT)>u},\\
&\Pi_3(u)=\pk{\sup_{t\in [\theta,1]}\inf_{s\in[te^{-2\delta T_u},t]}\LT(\sigma\int_0^{-\frac{1}{2\delta}\ln s} e^{-\delta v}dB(v)-\frac{c}{\delta}(1-{s}^{\frac{1}{2}})\RT)>u}.
\end{align*}
First we show the asymptotic of $\Pi_0(u)$.
For $u$ large enough
\begin{align*}
\Pi_0(u)&=\pk{\sup_{t\in [0, t_u+\lambda u^{-2}]}\inf_{s\in[te^{-2\delta T_u},t]}\LT(\sigma\int_0^{-\frac{1}{2\delta}\ln s} e^{-\delta v}dB(v)-\frac{c}{\delta}(1-{s}^{\frac{1}{2}})\RT)>u}\\
&=\pk{\sup_{t\in [0, t_u+\lambda u^{-2}]}\inf_{s\in[te^{-2\delta T_u},t]}\LT(\sigma\int_0^{-\frac{1}{2\delta}\ln s} e^{-\delta v}dB(v)-\frac{c}{\delta}(1-{s}^{\frac{1}{2}})\RT)>u}\\
&\leq\pk{\sup_{t\in [0, (1+\vn_1)(c^2/\delta^2+\lambda) u^{-2}]}\inf_{s\in[t(1+\vn_2)a,t]}\LT(\sigma\int_0^{-\frac{1}{2\delta}\ln s} e^{-\delta v}dB(v)-\frac{c}{\delta}(1-{s}^{\frac{1}{2}})\RT)>u}\\
&=\pk{\sup_{t\in [0, (1+\vn_1)(c^2/\delta^2+\lambda) ]}\inf_{s\in[(1+\vn_2)a,1]}\overline{Z}(st u^{-2})\frac{M_u(stu^{-2})}{M_u}>\frac{u}{M_u}}\\
&=:\Pi^{+\vn}_0(u),
\end{align*}
where $\overline{Z}(t)=\frac{Z(t)}{V_Z(t)}$, $a=e^{-2\delta T}$, $\vn_1\in(0,1)$ and $\vn_2\in\LT(0,(\frac{1}{a}-1)\wedge 1\RT)$ if $T\in(0,\IF]$, $\vn_2=0$ if $T=0$ .
Similarly,
\begin{align*}
\Pi_0(u)\geq\pk{\sup_{t\in [0, (1-\vn_1)(c^2/\delta^2+\lambda) ]}\inf_{s\in[(1-\vn_2)a,1]}\overline{Z}(st u^{-2})\frac{M_u(stu^{-2})}{M_u}>\frac{u}{M_u}}=:\Pi^{-\vn}_0(u).
\end{align*}
We have
\BQNY
1-\frac{M_u(t)}{M_u}=\frac{[G_u(t)V_Z(t_u)]^2-[G_u(t_u)V_Z(t)]^2}
{G_u(t)V_Z(t_u)[V_Z(t)G_u(t_u)+G_u(t)V_Z(t_u)]}.
\EQNY
and
\BQNY
[G_u(t)V_Z(t_u)]^2-[G_u(t_u)V_Z(t)]^2&=&\LT[\LT(u+\frac{c}{\delta}\RT)-\frac{c}{\delta}\sqrt{t}\RT]^2\frac{\sigma^2}{2\delta}(1-t_u)
-\LT[\LT(u+\frac{c}{\delta}\RT)-\frac{c}{\delta}\sqrt{t_u}\RT]^2\frac{\sigma^2}{2\delta}(1-t)\\
&=&\LT(u+\frac{c}{\delta}\RT)^2\frac{\sigma^2}{2\delta}(t-t_u)
-2\LT(u+\frac{c}{\delta}\RT)\frac{c\sigma^2}{2\delta^2}(\sqrt{t}-\sqrt{t_u})(1-t_u)-\frac{c^2\sigma^2}{2\delta^3}(t-t_u)\\
&=&\frac{\sigma^2}{2\delta}\LT[\LT(u+\frac{c}{\delta}\RT)^2-\LT(\frac{c}{\delta}\RT)^2\RT](\sqrt{t}-\sqrt{t_u})^2\\
&=&\frac{\sigma^2}{2\delta}\LT( u^2+\frac{2c}{\delta}u\RT)(\sqrt{t}-\sqrt{t_u})^2.
\EQNY
Since for any $t\in\Delta(u)$
\BQNY
\sqrt{\frac{\sigma^2}{2\delta}(1-t_u-\delta(u))}\leq V_Z(t)\leq\sqrt{\frac{\sigma^2}{2\delta}},\quad
u+\frac{c}{\delta}-\frac{c}{\delta}\sqrt{t_u+\delta(u)}\leq G_u(t)\leq u+\frac{c}{\delta},
\EQNY
then for all large $u$
\BQNY
V_Z(t_u)G_u(t)[G_u(t)V_Z(t_u)+V_Z(t)G_u(t_u)] \leq \frac{\sigma^2}{\delta}\LT(u+\frac{c}{\delta}\RT)^2
\EQNY
and
\BQNY
V_Z(t_u)G_u(t)[G_u(t)V_Z(t_u)+V_Z(t)G_u(t_u)]
&\geq&\frac{\sigma^2}{\delta}(1-t_u-\delta(u))\LT(u+\frac{c}{\delta}-\frac{c}{\delta}\sqrt{t_u+\delta(u)}\RT)^2\\
&\geq& \frac{\sigma^2}{\delta}\LT[\LT(u+\frac{c}{\delta}\RT)^2-u\RT].
\EQNY
Consequently, we have
\BQN\label{eq:MM12}
\lim_{u\rw\IF}\underset{s\in(0,1]}{\sup_{t\in(u^2\Delta(u))}}
\LT|\LT(1-\frac{M_u(stu^{-2})}{M_u}\RT)u^2-{\frac{1}{2}\LT(\sqrt{ts}-\frac{c}{\delta}\RT)^2 }\RT|=0.%\ \ \    \textbf{(Assumption A2)}
\EQN
For $0\leq t'\leq t<1$, the correlation function of $Z(t)$ equals
\BQN\label{rrr}
r(t,t')&=&\frac{\E{(\sigma\int_0^{-\frac{1}{2\delta}\ln t} e^{-\delta v}dB(v))(\sigma\int_0^{-\frac{1}{2\delta}\ln t'} e^{-\delta v}dB(v))}}{\sqrt{\frac{\sigma^2}{2\delta}(1-t)}\sqrt{\frac{\sigma^2}{2\delta}(1-t')}}\nonumber\\
&=&\frac{\sqrt{1-t}}{\sqrt{1-t'}}=1-\frac{t-t'}{\sqrt{1-t'}(\sqrt{1-t'}+\sqrt{1-t})},
\EQN
which implies that
\BQN\label{ar}
\sup_{t,t'\in\Delta(u),t'\neq t}\LT|\frac{1-r(t,t')}{\frac{1}{2}|t-t'|}-1\RT|
&=&\sup_{t,t'\in\Delta(u),t'\neq t}\LT|\frac{2}{\sqrt{1-t}(\sqrt{1-t'}+\sqrt{1-t})}-1\RT|\nonumber\\
&\leq&\frac{1}{1-(\frac{c}{c+\delta u})^2-(\frac{\ln u}{u})^2}-1\nonumber\\
&\rw& 0,\ \ u\rw\IF.
\EQN
For $t, t'\in\LT[0, (1+\vn_1)\LT(\frac{c^2}{\delta^2}+\lambda\RT)\RT]$ and $s, s'\in(0,1]$
\begin{align*}
&u^2\Var\LT(\overline{Z}(stu^{-2})\frac{M_u(stu^{-2})}{M_u}
-\overline{Z}(s't'u^{-2})\frac{M_u(s't'u^{-2})}{M_u}\RT)\\
&=\frac{u^2}{M^2_u}\E{\frac{Z(stu^{-2})}{1+\frac{c}{\delta u}(1-\sqrt{stu^{-2}})}
-\frac{Z(s't'u^{-2})}{1+\frac{c}{\delta u}(1-\sqrt{s't'u^{-2}})}}^2\\
&=\frac{u^2}{M^2_u}\E{\frac{Z(stu^{-2})-Z(s't'u^{-2})}{1+\frac{c}{\delta u}(1-\sqrt{stu^{-2}})}+\frac{Z(s't'u^{-2})}{1+\frac{c}{\delta u}(1-\sqrt{stu^{-2}})}
-\frac{Z(s't'u^{-2})}{1+\frac{c}{\delta u}(1-\sqrt{s't'u^{-2}})}}^2\\
&=\frac{u^2}{M^2_u}\LT(J_1(u)+J_2(u)+J_3(u)\RT),
\end{align*}
where
\BQNY
&&J_1(u)=\E{\LT(\frac{Z(stu^{-2})-Z(s't'u^{-2})}{1+\frac{c}{\delta u}(1-\sqrt{stu^{-2}})}\RT)^2},\\
&&J_2(u)=2\LT(\frac{1}{1+\frac{c}{\delta u}(1-\sqrt{stu^{-2}})}
-\frac{1}{1+\frac{c}{\delta u}(1-\sqrt{s't'u^{-2}})}\RT)\E{\frac{(Z(stu^{-2})-Z(s't'u^{-2}))Z(s't'u^{-2})}{1+\frac{c}{\delta u}(1-\sqrt{stu^{-2}})}}=0,\\
&&J_3(u)=\LT(\frac{1}{1+\frac{c}{\delta u}(1-\sqrt{stu^{-2}})}
-\frac{1}{1+\frac{c}{\delta u}(1-\sqrt{s't'u^{-2}})}\RT)^2\E{\LT(Z(s't'u^{-2})\RT)^2}.
\EQNY
Since for $t, t'\in\LT[0, (1+\vn_1)\LT(\frac{c^2}{\delta^2}+\lambda\RT)\RT]$ and $s, s'\in(0,1]$
\BQNY
&&\lim_{u\rw\IF}\frac{u^2}{M^2_u}J_1(u)=\lim_{u\rw\IF}\frac{u^2}{M^2_u(1+\frac{c}{\delta u}(1-\sqrt{stu^{-2}}))^2}\E{\LT(Z(stu^{-2})-Z(s't'u^{-2})\RT)^2}=\abs{st-s't'},\\
&&\lim_{u\rw\IF}\frac{u^2}{M^2_u}J_3(u)=\lim_{u\rw\IF}\frac{\sigma^2(1-s't'u^{-2})u^2}{2\delta M^2_u}\LT(\frac{\frac{c}{\delta u}\LT(\sqrt{s t u^{-2}}-\sqrt{s't' u^{-2}}\RT)}{\LT(1+\frac{c}{\delta u}(1-\sqrt{stu^{-2}})\RT)\LT(1+\frac{c}{\delta u}(1-\sqrt{s't'u^{-2}})\RT)}
\RT)^2=0,
\EQNY
we obtain
\BQN\label{VV}
\lim_{u\rw\IF}u^2\Var\LT(\overline{Z}(stu^{-2})\frac{M_u(stu^{-2})}{M_u}
-\overline{Z}(s't'u^{-2})\frac{M_u(s't'u^{-2})}{M_u}\RT)
&=&
|st-s't'|\nonumber\\
&=&2\Var\LT(\frac{1}{\sqrt{2}}B(st)-\frac{1}{\sqrt{2}}B(st)\RT).
\EQN
For some small $\theta\in(0,1)$, by \eqref{rrr} we obtain that for $t,t'\in[0,\theta]$
\BQN\label{hh}
\mathbb{E}\left(\overline{Z}(t)-\overline{Z}(t')\right)^2=2-2r(t,t')\leq\mathbb{C}_1|t-t'|%\ \  \ \textbf{(Assumption A5)}
\EQN
holds.
By \eqref{eq:MM12}, \eqref{ar},\eqref{VV}, \eqref{hh} and Lemma 5.1 in \cite{debicki2015parisian}, as $u\rw\IF$,
\begin{align*}
\Pi^{+\vn}_0(u)\sim \E{\sup_{t\in [0, (1+\vn_1)(c^2/\delta^2+\lambda) ]}\inf_{s\in[(1+\vn_2)a,1]}
\exp\LT(\frac{\sqrt{2\delta}}{\sigma}B(st)-\frac{\delta}{\sigma^2}st
-\frac{\delta}{\sigma^2}\LT(\sqrt{st}-\frac{c}{\delta}\RT)^2\RT)}\Psi\LT(\frac{u}{M_u}\RT),
\end{align*}
and
\begin{align*}
\Pi^{-\vn}_0(u)\sim \E{\sup_{t\in [0, (1-\vn_1)(c^2/\delta^2+\lambda) ]}\inf_{s\in[(1-\vn_2)a,1]}
\exp\LT(\frac{\sqrt{2\delta}}{\sigma}B(st)-\frac{\delta}{\sigma^2}st
-\frac{\delta}{\sigma^2}\LT(\sqrt{st}-\frac{c}{\delta}\RT)^2\RT)}\Psi\LT(\frac{u}{M_u}\RT).
\end{align*}
Letting $\vn_1,\vn_2\rw0$, we have
\BQN\label{PP0}
\Pi_0(u)\sim \E{\sup_{t\in [0, c^2/\delta^2+\lambda ]}\inf_{s\in[a,1]}
\exp\LT(\frac{\sqrt{2\delta}}{\sigma}B(st)-\frac{\delta}{\sigma^2}st
-\frac{\delta}{\sigma^2}\LT(\sqrt{st}-\frac{c}{\delta}\RT)^2\RT)}
\Psi\LT(\frac{u}{M_u}\RT), \ u\rw\IF.
\EQN

Next we show that
\begin{align*}
\Pi_1(u)=o\LT(\Pi_0(u)\RT),\quad \Pi_2(u)=o\LT(\Pi_0(u)\RT),\quad \text{and}\quad \Pi_3(u)=o\LT(\Pi_0(u)\RT).
\end{align*}
Let $Y(t), t\in \R$ be a stationary Gaussian process with continuous trajectories, unit variance and correlation function satisfying for a constant $\vn_3\in(0,\frac{1}{2})$
$$
r_{Y}(t)=1-\frac{(1+\vn_3)}{2}|t|.
$$
By \eqref{eq:MM12} and Slepian inequality in \cite{Pit96}, we have
\begin{align*}
\Pi_1(u)&\leq\sum_{i=1}^{N(u)}\pk{\sup_{t\in(t_u+ I_u(k))}\LT(\sigma\int_0^{-\frac{1}{2\delta}\ln t} e^{-\delta v}dB(v)-\frac{c}{\delta}(1-{t}^{\frac{1}{2}})\RT)>u}\\
&\leq\sum_{i=1}^{N(u)}\pk{\sup_{t\in (t_u+I_u(k))}\overline{Z}(t)>\mathcal{A}_u(k)}\\
&\leq\sum_{i=1}^{N(u)}\pk{\sup_{t\in (t_u+I_u(k))}Y(t)>\mathcal{A}_u(k)}\\
&=\sum_{i=1}^{N(u)}\pk{\sup_{t\in [0, \lambda]}Y(u^{-2}t)>\mathcal{A}_u(k)}
\end{align*}
where $\mathcal{A}_u(k):=\frac{u}{M_u}\LT(1+\frac{1-\vn_4}{2u^2}(\sqrt{u^2t_u+k\lambda }-c/\delta)^2-\frac{\vn_4}{u^2}\RT)$ and $\vn_4\in(0,1)$ is a small constant.
We observe that
\begin{align}\label{C1}
\inf_{1\leq k\leq N(u)}\mathcal{A}_u(k)\geq \frac{u}{M_u}\rw \IF, \ u\rw\IF.
\end{align}
Further,
\BQN\label{C2}
&&\lim_{u\rw\IF}\sup_{1\leq k\leq N(u)}\underset{t_1,t_2\in[0,\lambda]}{\sup_{t_1\neq t_2,}}
\LT|\mathcal{A}^2_u(k)\frac{\Var{\LT(Y(u^{-2}t_1)-Y(u^{-2}t_2)\RT)}}{\frac{2\delta(1+\vn_3)}{\sigma^2}|t_1-t_2|}-1\RT|\nonumber\\
&&\quad\quad=\lim_{u\rw\IF}\sup_{1\leq k\leq N(u)}\underset{t_1,t_2\in[0,\lambda]}{\sup_{t_1\neq t_2,}}
\LT|\mathcal{A}^2_u(k)\frac{2-2r_Y(u^{-2}t_1-u^{-2}t_2)}{\frac{2\delta(1+\vn_3)}{\sigma^2}|t_1-t_2|}-1\RT|\nonumber\\
&&\quad\quad=0,
\EQN
%where $a=\frac{2\delta^2 e^{-2\delta S}}{\sigma^2(1-e^{-2\delta S})^2}$,
and
\BQN\label{C3}
&&\sup_{1\leq k\leq N(u)}\underset{t_1,t_2\in[0,\lambda]}{\sup_{|t_1-t_2|<\epsilon}}
\mathcal{A}^2_u(k)\E{\LT(Y(u^{-2}t_1)-Y(u^{-2}t_2)\RT)Y(0)}\nonumber\\
&&\quad\quad\leq \mathbb{C}_2 u^2\underset{t_1,t_2\in[0,\lambda]}{\sup_{|t_1-t_2|<\epsilon}}\LT|r_Y(u^{-2}t_1)-r_Y(u^{-2}t_2)\RT|\nonumber\\
&&\quad\quad\leq \mathbb{C}_3 u^2\underset{t_1,t_2\in[0,\lambda]}{\sup_{|t_1-t_2|<\epsilon}}\LT|\frac{1+\vn_3}{2}u^{-2}(t_1-t_2)\RT|\nonumber\\
&&\quad\quad\leq \mathbb{C}_4 \underset{t_1,t_2\in[0,\lambda]}{\sup_{|t_1-t_2|<\epsilon}}|t_1-t_2|\rw 0 ,\ \ u\rw\IF, \ \epsilon\rw0.
\EQN
According to \eqref{C1}, \eqref{C2}, \eqref{C3} and Lemma 5.3 of \cite{KEP2015},  we have
as $u\rw\IF, \vn_4\rw 0,\ \lambda\rw\IF$
\BQN\label{bound3}
\Pi_1(u)&\leq& \mathbb{C}_5\lambda\sum_{k=1}^{N(u)}\Psi\LT(\mathcal{A}_u(k)\RT)\nonumber\\
&\sim& \mathbb{C}_5\lambda\sum_{k=1}^{N(u)}\frac{1}{\sqrt{2\pi}\mathcal{A}_u(k) }e^{-\frac{\mathcal{A}^2_u(k)}{2}}\nonumber\\
&\leq& \mathbb{C}_5\lambda\sum_{k=1}^{N(u)}\frac{M_u}{\sqrt{2\pi}u }
\exp\LT(-\frac{u^2}{2M_u^2}\LT(1+\frac{1-\vn_4}{u^2}\LT(\sqrt{u^2t_u+k\lambda }-c/\delta\RT)^2-\frac{2\vn_4}{u^2}\RT)\RT)\nonumber\\
&\sim& \mathbb{C}_5\lambda \Psi\LT(\frac{u}{M_u}\RT)e^{\frac{\vn_4}{M_u^2}}\sum_{k=1}^{N(u)}
\exp\LT(-\frac{1-\vn_4}{2M_u^2}\LT(\sqrt{u^2t_u+k\lambda }-c/\delta\RT)^2\RT)\nonumber\\
&\leq&\mathbb{C}_6 \Psi\LT(\frac{u}{M_u}\RT)e^{\frac{\sigma^2\vn_4}{2\delta}} \lambda\sum_{k=1}^{\IF}e^{-\mathbb{C}_{7} k\lambda}= o\LT(\Psi\LT(\frac{u}{M_u}\RT)\RT).
\EQN

Moreover, for all $u$ large
\BQN\label{assA41}
\frac{1}{M_u(t)}-\frac{1}{M_u}&\geq& \frac{[G_u(t)V_Z(t_u)]^2-[G_u(t_u)V_Z(t)]^2}{2uV_Z^3(t_u)G_u(t_u)}\nonumber\\
&=&\frac{\frac{\sigma^2}{2\delta}(u^2+\frac{2c}{\delta}u)(\sqrt{t}-\sqrt{t_u})^2}{2u[\frac{\sigma^2}{2\delta}(1-t_u)]^{3/2}[u+\frac{c}{\delta}(1-\sqrt{t_u})]}\nonumber\\
%&=&\frac{\sqrt{2}(\delta u+c)^4}{2u\delta^2\sigma(\delta u^2+2cu)^{3/2}}\frac{t^2}{(\sqrt{t+t_u}+\sqrt{t_u})^2}\nonumber\\
&\geq&\mathbb{C}_8 (\sqrt{t}-\sqrt{t_u})^2\nonumber\\
&\geq&\frac{\mathbb{C}_8\LT(\frac{\ln u}{u}\RT)^4}{\LT(\sqrt{\LT(\frac{\ln u}{u}\RT)^2+(\frac{c}{\delta u+c})^2}+\frac{c}{\delta u+c}\RT)^2}\nonumber\\
&\geq&\mathbb{C}_8\frac{(\ln u)^{2}}{u^2}%\ \ \textbf{(Assumption A4)}
\EQN
holds for any $t\in\LT[t_u+\delta(u),\theta\RT]$, therefore
\BQNY
\sup_{t\in\LT[t_u+\delta(u),\theta\RT]}M_u(t)\leq \LT(\frac{1}{M_u}+\mathbb{C}_8\frac{(\ln u)^2}{u^2}\RT)^{-1}.
\EQNY
\COM{Further, by \eqref{hh} we have
\begin{align*}
\E{(\overline{Z}(t_1)-\overline{Z}(t_2))^2}\leq \mathbb{C}_9 \abs{t_1-t_2}, t_1,t_2\in[0,\theta].
\end{align*}}
Thus the above inequality combined with \eqref{hh} and Theorem 8.1 in \cite{Pit96} derives that
\BQN\label{bound4}
\Pi_2(u)&\leq&\pk{\sup_{t\in [t_u+\delta(u),\theta]}\overline{Z}(t )M_u(t)>u}\nonumber\\
&\leq&\pk{\sup_{t\in [0,\theta]}\overline{Z}(t )>u\LT(\frac{1}{M_u}+\mathbb{C}_8\frac{(\ln u)^{2}}{u^2}\RT)}\nonumber\\
&\leq&\mathbb{C}_{9} u^2\Psi\LT(u\LT(\frac{1}{M_u}+\mathbb{C}_8\frac{(\ln u)^{2}}{u^2}\RT)\RT)\nonumber\\
&\leq& o\LT(\Psi\LT(\frac{u}{M_u}\RT)\RT), \ u\rw\IF.
\EQN

Finally, since
\BQNY
\sup_{t\in [\theta,1]}V^2_Z(t)
%=\sup_{t\in [\theta,1]}\frac{\sigma^2u^2(1-t)}{2\delta M_u^2( u+\frac{c}{\delta}(1-t^{1/2}))^2}
\leq \frac{\sigma^2}{2\delta}(1-\theta),\ \ \hbox{and} \ \
\E{\sup_{t\in [\theta,1]}Z(t)}\leq \mathbb{C}_{10}<\IF,
\EQNY
by Borell inequality in \cite{AdlerTaylor}
\BQN
\Pi_3(u)\leq\pk{\sup_{t\in [\theta,1]}Z(t)>u}
\leq \exp\LT(-\frac{\delta(u-\mathbb{C}_{10})^2}{\sigma^2(1-\theta)}\RT)
%=o\LT(\Psi\LT(\frac{1}{\sigma}\sqrt{2\delta u^2+4cu}\RT)\RT)
=o\LT(\Psi\LT(\frac{u}{M_u}\RT)\RT),
\ u\rw \IF,
\EQN
which combined with \eqref{bound1}, \eqref{bound2}, \eqref{PP0}, \eqref{bound3} and \eqref{bound4} shows that
$$\psi(u)\sim \Pi_0(u),\quad u\rw\IF.$$
Consequently, letting $\lambda\rw\IF$, we have
\begin{align*}
\psi(u)&\sim \E{\sup_{t\in [0, \IF)}\inf_{s\in[a,1]}
\exp\LT(\frac{\sqrt{2\delta}}{\sigma}B(st)-\frac{\delta}{\sigma^2}st
-\frac{\delta}{\sigma^2}\LT(\sqrt{st}-\frac{c}{\delta}\RT)^2\RT)}
\Psi\LT(\frac{1}{\sigma}\sqrt{2\delta u^2+4cu}\RT)\\
&= \E{\sup_{t\in [0, \IF)}\inf_{s\in[a,1]}
\exp\LT(\sqrt{2}B(st)-st
-\LT(\sqrt{st}-\frac{c}{\sigma\sqrt{\delta}}\RT)^2\RT)}
\Psi\LT(\frac{1}{\sigma}\sqrt{2\delta u^2+4cu}\RT), \ u\rw\IF.
\end{align*}
\QED

\prooftheo{ruintime}
We use the same notation as in the proof of \netheo{paripan1}. For $x\in\LT(-\frac{c^2}{\delta^2},\IF\RT)$ and $u>0$
\begin{align*}
&\pk{u^{2}\LT(e^{-2\delta\eta(u)}-\LT(\frac{c}{\delta u+c}\RT)^2\RT)\leq x\big| \eta(u)<\IF}\\
&=\frac{\pk{\inf_{t\in[-\frac{1}{2\delta}\ln \LT(t_u+u^{-2}x\RT),\IF)}\sup_{s\in[t,t+T_u]}\widetilde{R}_u^{\delta}(s)<0}}
{\pk{\inf_{t\in[0,\IF)}\sup_{s\in[t,t+T_u]}\widetilde{R}_u^{\delta}(s)<0}}\\
&=\frac{\mathbb{P}\left\{\sup_{t^*\in[0,t_u+u^{-2}x]}\inf_{s^*\in[t^*e^{-2\delta T_u},t^*]}\LT(\sigma\int_0^{-\frac{1}{2\delta}\ln s^*} e^{-\delta v}dB(v)-c\int_0^{-\frac{1}{2\delta}\ln s^*} e^{-\delta v}dv\RT)>u\right\}}
{\mathbb{P}\left\{\sup_{t^*\in[0,1]}\inf_{s^*\in[t^*e^{-2\delta T_u},t^*]}\LT(\sigma\int_0^{-\frac{1}{2\delta}\ln s^*} e^{-\delta v}dB(v)-c\int_0^{-\frac{1}{2\delta}\ln s^*} e^{-\delta v}dv\RT)>u\right\}}\\
&=\pk{u^{2}\LT(\tau^*_u-t_u\RT)\leq x\big| \tau^*_u<1},
\end{align*}
where $\tau^*_u=\{t\geq 0: \sigma\int_0^{-\frac{1}{2\delta}\ln t^*} e^{-\delta v}dB(v)-\frac{c}{\delta}(1-{t^*}^{\frac{1}{2}})>u\}$.\\
For $\psi_x(u):=\mathbb{P}\left\{\sup_{t^*\in[0,t_u+u^{-2}x]}\inf_{s^*\in[t^*e^{-2\delta T_u},t^*]}\LT(\sigma\int_0^{-\frac{1}{2\delta}\ln s^*} e^{-\delta v}dB(v)-c\int_0^{-\frac{1}{2\delta}\ln s^*} e^{-\delta v}dv\RT)>u\right\}$, using the similar argumentation about $\Pi_0(u)$ as in the proof of \netheo{paripan1} with $\lambda=x$, we obtain
\begin{align*}
\psi_x(u)&\sim\E{\sup_{t\in [0, c^2/\delta^2+x ]}\inf_{s\in[a,1]}
\exp\LT(\frac{\sqrt{2\delta}}{\sigma}B(st)-\frac{\delta}{\sigma^2}st
-\frac{\delta}{\sigma^2}\LT(\sqrt{st}-\frac{c}{\delta}\RT)^2\RT)}
\Psi\LT(\frac{1}{\sigma}\sqrt{2\delta u^2+4cu}\RT)\\
&=\E{\sup_{t\in [0,\frac{c^2}{\sigma^2\delta}+\frac{\delta x}{\sigma^2} ]}\inf_{s\in[a,1]}
\exp\LT(\sqrt{2}B(st)-st
-\LT(\sqrt{st}-\frac{c}{\sigma\sqrt{\delta}}\RT)^2\RT)}
\Psi\LT(\frac{1}{\sigma}\sqrt{2\delta u^2+4cu}\RT), \ u\rw\IF.
\end{align*}
Thus
\begin{align*}
\pk{u^{2}\LT(e^{-2\delta\eta_u}-\LT(\frac{c}{\delta u+c}\RT)^2\RT)\leq x\big| \eta_u<\IF}
=\frac{\psi_x(u)}{\psi(u)}
\sim\frac{\PP^f_a[0,\frac{c^2}{\sigma^2\delta}+\frac{\delta x}{\sigma^2}]}{\PP^f_a[0,\IF)}, \ u\rw\IF.
\end{align*}
\QED

{\bf Acknowledgement}: Thanks to  Swiss National Science Foundation Grant no.  200021-166274.
\bibliographystyle{plain}
\bibliography{InfiniteParisian}
\end{document}